\documentclass[10pt]{article}
\usepackage{amssymb}
\usepackage{amsmath}
\usepackage{epsfig}

\newtheorem{theorem}{Theorem}[section]

\newtheorem{definition}[theorem]{Definition}

\newtheorem{remark}[theorem]{Remark}
\newtheorem{example}[theorem]{Example}

\newcommand{\boproof}{\noindent {\bf Proof. }}
\newcommand{\eoproof}{\hspace*{\fill} $\square$ \vspace{5pt}}

\newcommand{\R}{\mathbb R}

\newcommand{\Z}{\mathbb Z}

\newcommand{\Orthant}{\mathbb O}

\DeclareMathOperator{\cone}{cone}

\DeclareMathOperator{\conv}{conv}

\DeclareMathOperator{\maxdeg}{maxdeg}

\begin{document}
\title{Integral Function Bases}
\author{Raymond Hemmecke, Robert Weismantel\\Otto-von-Guericke-University Magdeburg}
\date{}
\maketitle

\begin{abstract}
Integral bases, a minimal set of solutions to $Ax\leq b, x\in\Z^n$
that generate any  other solution to $Ax\leq b, x\in\Z^n$, as a
nonnegative integer linear combination, are always finite and are
at the core of the Integral Basis Method introduced by Haus,
K{\"o}ppe and Weismantel.

In this paper we present one generalization of the notion of
integral bases to the nonlinear situation with the intention of
creating an integral basis method also for nonlinear integer
programming.
\end{abstract}


\section{Introduction}
In the past fifty years many efforts have been undertaken to study
linear integer optimization problems from different mathematical
and algorithmic viewpoints. As a result, a basic understanding of
the geometry of integer programming problems defined by linear
equations and/or linear inequalities is present today. This
knowledge has been partly turned into algorithmic tools to tackle
discrete optimization problems in practice.

The attempts to study the geometry of integer points in polyhedral
sets are based on two basic mathematical concepts. One is the
notion of a lattice. More precisely,  a basis of a lattice
 $L$ is a subset of linearly independent vectors that
allows one to  generate all points in the lattice with respect to
taking integer linear combinations. The geometric properties of
particular bases in lattices made it possible to design algorithms
for solving specific linear integer programming  problems, mainly
problems without lower and upper bounds on the variables  and
linear problems with a fixed number of discrete variables
\cite{lenstra,Schrijver:86}. The notion of bases of a lattice can
be further refined so as to yield so-called integral generating
sets for cones and polyhedra. Roughly speaking, integral
generating sets  extend -- besides lattices -- the notions of
extreme points and rays in polyhedra and cones to integer  points
in such sets. More precisely, an integral generating set for a set
$S \subseteq \Z^n$ is a subset of $S$ with the property that every
member of $S$ can be represented as a nonnegative integer
combination of the elements in $S$. Of course, $S$ itself
constitutes an integral generating set of itself. The key question
is to detect an integral generating set that is finite and minimal
with respect to inclusion. This immediately raises the question to
characterize those sets $S$ of lattice points that possess a
finite integral generating set. This question is answered in
Section \ref{Section: Integral Bases} of this paper.

Indeed, integral generating sets  have important implications for
the theory of linear integer programming. Most importantly,
optimality conditions for integer  optimization problems can be
derived through integral generating sets. Such sets also provide a
basic understanding of integral polyhedra and  totally dual
integral systems of inequalities \cite{Giles+Pulleyblank:79}. Last
but not least they play a central role in the development of
integer simplex type methods of linear integer programs
\cite{Haus+Koeppe+Weismantel}. In fact, it is quite obvious to see
that if a finite integral generating set for a discrete set of
points is available, then we can reformulate the problem of
detecting a particular element in $S$ as the problem of  detecting
a nonnegative integer multiplier associated with the new
representation through an integral generating set. Integral
generating sets therefore allow a new representation of the same
set $S$ in some other space. The beautiful fact is that if we
start off with a set $S$ that is the feasible region of an integer
linear program in nonnegative variables, then also after
reformulation the new optimization problem happens to be a linear
integer program in nonnegative variables. This follows simply from
the fact that the integral generating set enables us to express
every point as a nonnegative integer combination.

Suppose now that $S$ does not have a finite integral generating
set. This in fact may happen even though $S$ corresponds to all
the integer points in a polyhedron.  In particular, if the
constraints defining $S$ are not linear, even in quite restrictive
cases  $S$ does not possess a finite integral generating set.

Then the idea to use integral generating sets for reformulation
issues is not possible, because the generating set is infinite. We
can simply not write down any finite representation of the
reformulated problem. In order to cope with this scenario, it
requires to generalize the notion of an integral generating set
from the linear case to a nonlinear setting. We  will refer to
such sets as integral function bases, since they enable us to
derive representations by means of  nonnegative polynomial
combinations  instead of nonnegative linear combinations. This is
the central topic of Section \ref{Section: Nonlinear Integral
Bases. Definition and Motivation}. In turn, our generalization
allows us to formulate optimality conditions for integer
polynomial programming problems.

We also analyze the situation when an integral function basis for
the integer points in  a complicated semi-algebraic set is
replaced by the condition of being members of a relaxation of the
semi-algebraic set itself. This question is in particular
motivated by the design of pivoting type methods for polynomial
integer and mixed integer programming, a topic that we regard as
theoretically and practically challenging, but important.

\section{Integral Bases}\label{Section: Integral Bases}

Let us start by defining the notion of an integral basis.

\begin{definition}\label{Integral Basis}
Let $S\subseteq\Z^n$. Then we call $T\subseteq S$ an
{\bf integral basis} of $S$, if for every $s\in S$ there exists a
finite (integer) linear combination $s=\sum \alpha_i t_i$ with
$t_i\in T$ and $\alpha_i\in\Z_+$.
\end{definition}

\begin{figure}[tbh]
\begin{center}
\epsfig{file=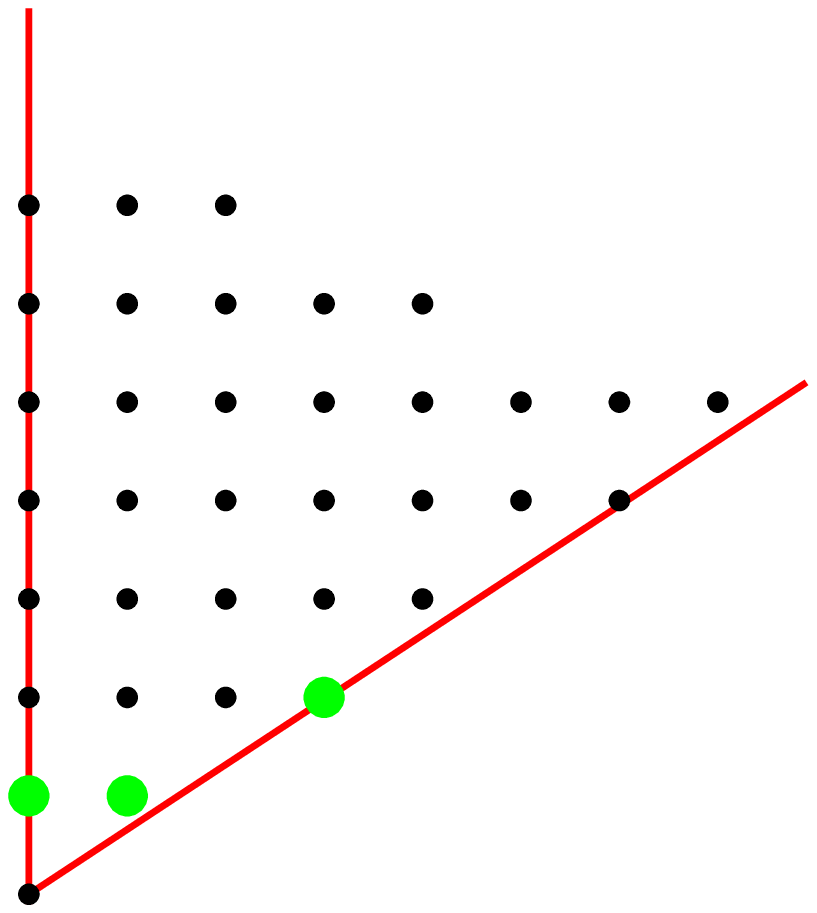, height=3cm}\hspace{3cm}
\epsfig{file=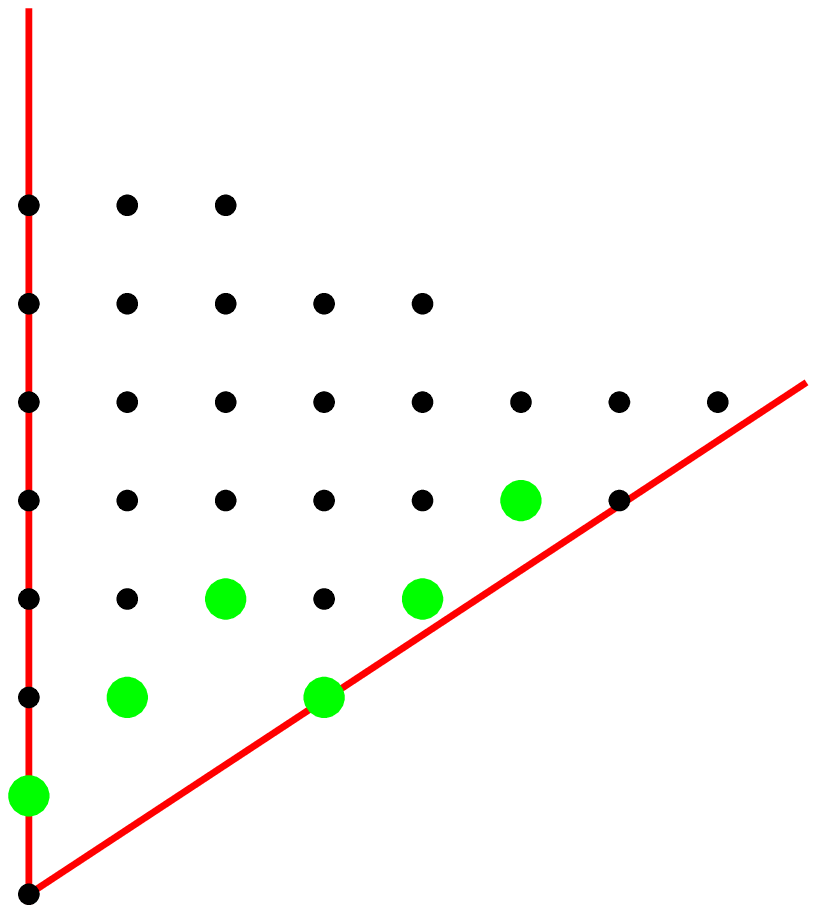, height=3cm}\\[.3cm]
\refstepcounter{figure} \label{Two integral bases} Figure \thefigure:
Minimal integral bases of two sets of lattice points
\end{center}
\vspace*{-0.3cm}
\end{figure}

Note that an integral basis of $S$ is allowed to contain elements only
from $S$ itself! With this definition, Bertsimas and Weismantel
\cite{Bertsimas+Weismantel} showed the following characterization of
which rational polyhedra (or more precisely the integer points in such
polyhedra) have a  {\it finite} integral basis.

\begin{theorem}[Bertsimas and Weismantel \cite{Bertsimas+Weismantel}]
\label{Finite integral basis for Ax<=b} For $A\in\Z^{d\times n}$
and $b\in\Z^d$, define the sets $P=\left\{x\in\R_+^n: Ax\leq
b\right\}$, $S=P\cap \Z^n$, and $C=\{x\in\R_+^n: Ax\leq 0\}$.
\begin{itemize}
\item[{\bf (a)}]
There exists a finite integral generating set of $S$ if and only
if $S$ contains all  but finitely many integer points in $C\cap
\Z^n_+$.
\item[{\bf (b)}]
If a finite integral generating set of $S$ exists, then there is a
unique integral basis of $S$.
\end{itemize}
\end{theorem}

Note that for $b=0$, this theorem simply states existence (and
uniqueness) of a Hilbert basis for the (pointed) rational polyhedral
cone $\{x\in\R_+^n: Ax\leq 0\}$.

Now let us give a novel and more general characterization of which
sets of lattice points have a finite integral basis. As we do not make
any structural assumption on the set of lattice points, we have to be
cautious to check whether the integral bases that we construct do
indeed consist of lattice points from our original sets only.

\begin{theorem}\label{General Finiteness Theorem}
Let $S\subseteq\Z^n$ be any set of lattice points in $\Z^n$.
\begin{itemize}
\item[{\bf (a)}]
$S$ has a finite integral basis if and only if $C=\cone(S)$ is a
  rational polyhedral cone.
\item[{\bf (b)}]
If the cone $C=\cone(S)$ is rational and pointed, there is a unique
finite integral basis that is minimal with respect to set inclusion.
\end{itemize}
\end{theorem}

\boproof Let us start showing part (a). If $C=\cone(S)$ is not a
rational polyhedral cone, $S$ cannot have a finite integral basis
$G\subseteq S$, since $C=\cone(S)=\cone(G)$ would be a rational
cone, contradicting our initial assumption on $C$.

Now we show the remaining claim that $S$ has a finite integral basis
if $C=\cone(S)$ is rational by explicitly constructing such a finite
basis. It should be noted that this integral basis need not be
minimal.

First, let us triangulate $C$ into (finitely many!) simplicial
cones $C_1,\ldots,C_k$. Note that we can and do choose such a
triangulation for which the generators of the cones $C_i$ are also
among the (finitely many) generators of $C$. Thus, as $C$, each
cone $C_i$ is generated by (finitely many) elements $S_i$ of $S$.
It remains to show that for each rational simplicial cone
$C_i=\cone(S_i)$, the set $C_i\cap S$ has a finite integral basis
$G_i$. Then the union of all $G_i$, $i=1,\ldots,k$, is clearly a
finite integral basis for $S$.

For $C_i=\cone(S_i)$ and $S_i=\{v_1,\ldots,v_r\}$, consider the
parallelepiped
\[
F=\left\{\sum_{j=1}^r \alpha_j v_j: 0\leq\alpha_1,\ldots,\alpha_r<1\right\}.
\]
As $F_i$ is bounded, $F$ contains only finitely many lattice
points $\{f_1,\ldots,f_t\}$ in $\Z^n$. Moreover, $C_i\cap\Z^n$ is
the disjoint union of the following $t$ sets $F_1,\ldots,F_t$ with
\[
F_j=\left\{f_j+\sum_{j=1}^r \alpha_j v_j: \alpha_1,\ldots,\alpha_r\in\Z_+\right\}.
\]
We construct now a finite integral basis for $C_i\cap S$.

Consider any $F_j$, $j=1,\ldots,t$. As $C_i$ is a simplicial cone,
each point in $F_j$ has a unique representation as
$f_j+\sum_{j=1}^r\alpha_j v_j$ implying that there is a one-to-one
correspondences $\phi_j$ between $F_j$ and $\Z_+^r$ given by
\[
\phi_j\left(f_j+\sum_{j=1}^r \alpha_j v_j\right)=(\alpha_1,\ldots,\alpha_r).
\]
To construct a finite integral basis for $F_j\cap S$, consider the
set $\phi_j(F_j\cap S)\subseteq\Z_+^r$. By the Gordan-Dickson
Lemma, there are only finitely many points $\{g_1,\ldots,g_p\}$
that are minimal with respect to the partial ordering $\leq$
defined on $\Z_+^r$. Thus, each point $\phi_j(F_j\cap S)$ can be
written as a positive integer linear combination of
$g_1,\ldots,g_p$ and of the unit vectors $e_1,\ldots,e_r$. Thus,
every element in $F_j\cap S$ is a positive integer linear
combination of $\phi^{-1}(g_1),\ldots,\phi^{-1}(g_p)\in F_j\cap S$
together with $\phi^{-1}(e_1),\ldots,\phi^{-1}(e_r)\in
S_i\subseteq S$. Let $G_{i,j}$ denote the set of all these
vectors. Clearly, the union $G_i$ over all $G_{i,j}$, $j=1,\ldots,
t$, forms a finite integral basis for $C_i\cap S$, and claim (a)
is proved.

Let us prove claim (b) now. As $\cone(S)$ is pointed, there is some
vector $c\in\R^n$ such that
$\{x\in\R^n: c^\intercal x\leq 0\}\cap\cone(S)=\{0\}$. Assume that
$U=\{u_1,\ldots,u_r\}$ and $V=\{v_1,\ldots,v_t\}$ are two different
inclusion minimal integral bases of $S$. Moreover, assume that
w.l.o.g. $u_1\not\in V$. Minimality of $U$ implies that $u_1$ cannot
be written as a positive integer linear combination of elements in
$U\setminus\{u_1\}$. However, as $V$ is an integral basis of $S$ and
$u_1\in S$, there is a nonnegative integer linear combination
$u_1=\sum_{j=1}^t\alpha_j v_j$. Clearly, as $u_1\not\in V$ and as the
coefficients are nonnegative integers, we have
$c^\intercal v_j<c^\intercal u_1$ whenever $\alpha_j>0$. As also $U$
is an integral basis of $S$ and as all $v_j\in S$, there are
nonnegative integer linear combinations $v_j=\sum_{i=1}^r\beta_{i,j}
u_i$. Moreover, $c^\intercal u_i\leq c^\intercal v_j$ whenever
$\beta_{i,j}>0$. Plugging these representations into
$u_1=\sum_{j=1}^t\alpha_j v_j$, we get a representation of $u_1$ as a
nonnegative integer linear combination of elements in $U$. However, by
construction, they all have a scalar product with $c$ that is strictly
less than $c^\intercal u_1$. Thus, we have written $u_1$ as a
nonnegative integer linear combination of elements in
$U\setminus\{u_1\}$, a contradiction to our assumption that $U$ is a
set inclusion minimal integral basis, and the claim is proved.
\eoproof

Note that for sets of the form $\{x\in\Z^n: Ax\leq 0\}$,
Theorem~\ref{General Finiteness Theorem} again simply states
existence of finite Hilbert bases for rational polyhedral cones
and uniqueness of the minimal Hilbert basis if the cone is
pointed. It is easy to show that the minimal Hilbert basis of a
cone must consist of lattice points from the fundamental
parallelepiped and is thus finite. The tricky part for the proof
of Theorem~\ref{General Finiteness Theorem} was the fact, that not
all points of this parallelepiped could be assumed to belong to
$S$. The two examples in Figure~\ref{Two integral bases},
page~\pageref{Two integral bases} already illustrate this
difficulty.

Let us now show how Theorem~\ref{General Finiteness Theorem} implies
the special case, Theorem~\ref{Finite integral basis for Ax<=b}.

{\bf Proof} of Theorem~\ref{Finite integral basis for Ax<=b}. Let
us show part (a) first. If $S$ is finite, nothing is left to show.
Thus, assume that $S$ is not finite and therefore also
$C\neq\{0\}$. Assume that $S$ contains all but finitely many
integer points in $C\cap\Z^n_+$. In particular, $S$ contains an
(integer) point of every extreme ray of $C$. By Minkowski's
theorem, we have $\conv(P\cap\Z^n)=\conv(G)+C$, where $G\subseteq
S$ is the set of extreme points in $\conv(P\cap\Z^n)\subseteq P$.
(Since $P$ does not contain a line, $G\neq\emptyset$.) Thus,
$\cone(S)=\conv(G)+C$, as $G\subseteq S$ and as $S$ contains an
(integer) point of every extreme ray of $C$. Consequently,
$\cone(S)$ is a rational cone and thus has a finite integral basis
by Theorem~\ref{General Finiteness Theorem}.

Now assume that there are infinitely many integer points in $C$
that do not belong to $S$. In particular, $0\not\in P$ as
otherwise $0=A0\leq b$ implying $C\subseteq P$ and thus
$C\cap\Z^n\subseteq S$. Assume for the moment that each extreme
ray of $C$ contains a (nonzero!) point of $P$. Fix any extreme ray
of $C$ and let $v\in C\cap P$ be a point on this ray. Then any
point $w=kv$, $k\geq 1$, on this ray must belong to $P$. This
follows from $Aw=k(Av)=(k-1)v+v\leq (k-1)\cdot 0+b=b$ and $Av\leq
0$, $Av\leq b$, and $k\geq 1$. Therefore, $w\in P$ as claimed. By
convexity of $P$, $P$ must contain the convex hull $H$ of all
these half-lines $\{kv:k\geq 1, k\in\R\}$. As $C\setminus H$ is
bounded, only a finite number of integer points in $C$ can lie in
$C\setminus H$. As $S=P\cap\Z^n$, this implies that only finitely
many integer points $C$ can lie outside of $S$, contradicting our
initial assumption on $C$. This implies that there must be an
extreme ray $R$ of $C$ that does not contain any point of $P$.

We now show that $\cone(S)$ cannot be a rational cone, and the
result follows again by Theorem~\ref{General Finiteness Theorem}.
Assume on the contrary that $\cone(S)$ is a rational cone. By
convexity of $S$, every ray in $\cone(S)$ has a nontrivial
intersection with $S$ and thus also with $P$. This implies that
the extreme ray $R$ does not belong to $\cone(S)$. As $\cone(S)$
is rational, there exists a finite (rational) description
\[
\cone(S)=\{x\in\R^n: c_i^\intercal x\leq 0,i=1,\ldots,p\}.
\]
Let $v$ be any rational vector with $R=\cone(v)$. Then
$v\not\in\cone(S)$ implies that there is some index $j$ such that
$c_j^\intercal v>0$. Now consider any integer point $w\in S$. As
$S=\conv(G)+C$, all integer points on the half-line $\{w+\alpha
v:\alpha\geq 0\}$ belong to $S$. Moreover, as $v$ is a rational
vector, there are infinitely many integer points on this
half-line. However, as $c_j^\intercal v>0$, we have $c_j^\intercal
(w+\alpha v)>0$ for sufficiently large $\alpha$, implying that
there are integer points of $S$ that lie outside of $\cone(S)$.
This contradiction shows that $\cone(S)$ is not a rational cone
and part (a) is proved.

As part (b) of our claim follows now immediately from part (b) of
Theorem~\ref{General Finiteness Theorem}, nothing is left to
show. \eoproof

A natural question that we may ask is, whether there are other
special cases of interesting sets of lattice points that have a
finite integral basis by Theorem~\ref{General Finiteness Theorem}.

One natural guess would be the integral points in a convex region.
However, convexity alone is not enough to ensure that $\cone(S)$
is rational, as can be seen by looking at a polyhedral cone with
irrational generators. Thus, some notion of ``rational
generators'' of the region should be defined.

\begin{figure}[tbh]\label{cone(S) is not rational}
\begin{center}
\epsfig{file=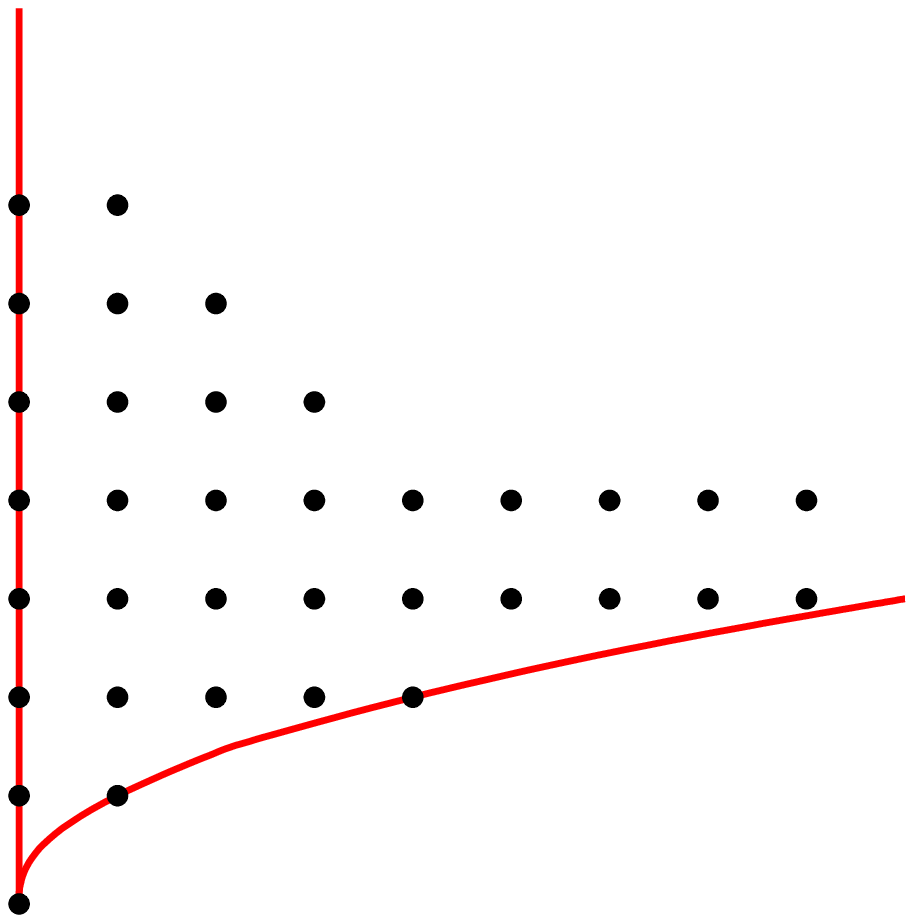, height=3cm}\\[.3cm]
\refstepcounter{figure} Figure \thefigure:
$\cone(S)$ is not rational.
\end{center}
\vspace*{-0.3cm}
\end{figure}

With this in mind, we may try to look at sets that are parametrized by
convex polynomials that have rational coefficients only. Again, there
is a simple counter-example. For the lattice points $S$ in the
parametrized set
\[
\left\{ \left(\begin{array}{c}
x\\
y\\
\end{array}
\right)
\in\R^2:
\left(\begin{array}{c}
x\\
y\\
\end{array}
\right)
=
\left(\begin{array}{c}
s\\
s^2\\
\end{array}
\right)
+\left(\begin{array}{c}
0\\
t\\
\end{array}
\right), s,t\in\R_+ \right\},
\]
we easily see that $\cone(S)$ is not rational, see Figure
\ref{cone(S) is not rational}.

We conclude that even convexity {\bf and} rationality of
generators is generally not enough to ensure finiteness of an
integral basis. It can be shown that under the assumption that if
in addition the given set itself is convex and that if it contains
the unit vectors of the positive orthant, a finite integral basis
does exist.

\section{Nonlinear Integral Bases. Definition and Motivation}
\label{Section: Nonlinear Integral Bases. Definition and
Motivation}

Theorem \ref{General Finiteness Theorem} characterizes when linear
integral bases exist. What can we do if the conditions of the
theorem do not hold? For instance, if we consider the set
$S=\{(x,y)\in\Z^2_+:y\geq 1\}$.

\begin{figure}[tbh]
\begin{center}
\epsfig{file=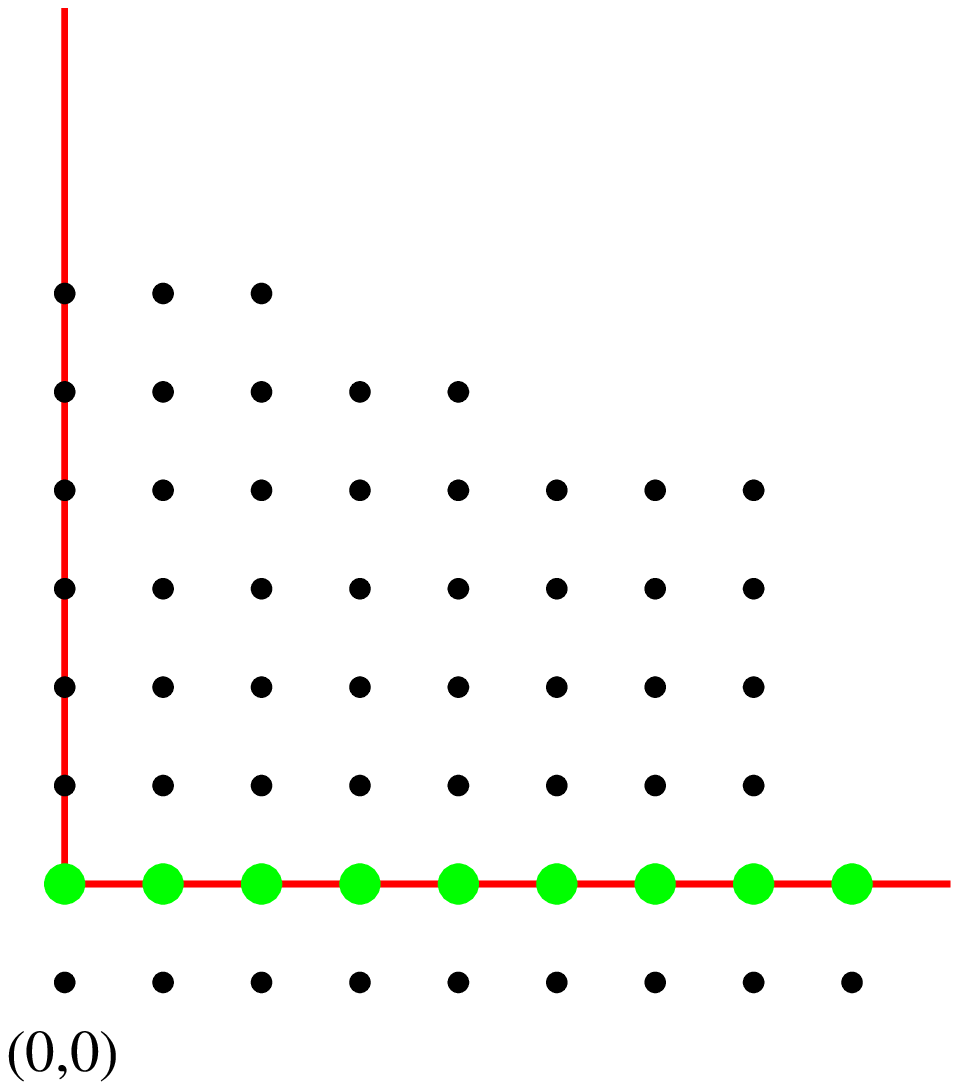, height=3cm}\\[.3cm]
\refstepcounter{figure} \label{IntegralBasisNotFinite} Figure
\thefigure: Example of an infinite integral basis
\end{center}
\vspace*{-0.3cm}
\end{figure}

In this case, the set $\cone(S)$ is not finitely generated, and
thus there does not exist a finite integral basis of $S$. For
obtaining a finite representation in this example, it becomes
necessary to extend the notion of an integral basis to -- what we
call -- an integral function basis. Our goal then becomes to
identify sets of points that have a finite integral function
basis. In the following we consider sets
$S=\{y\in\R^n:y=g(\lambda),\lambda\in\Z^d_+\}\cap\Z^n$, where
$g:\R^d\rightarrow\R^n$ is a vector of functions with components
$g_i:\R^d\rightarrow\R$, $i=1,\ldots,n$, and with
$g(\Z^d)\subseteq\Z^n$.

Note that when all $g_i$ are linear functions, our set $S$
corresponds to the lattice points of a rational polyhedral cone.
Other possible functions are polynomials in $\Z[\lambda]$, certain
stair-case functions, or even suitable combinations of all $3$
types.

\begin{figure}[tbh]
\begin{center}
\epsfig{file=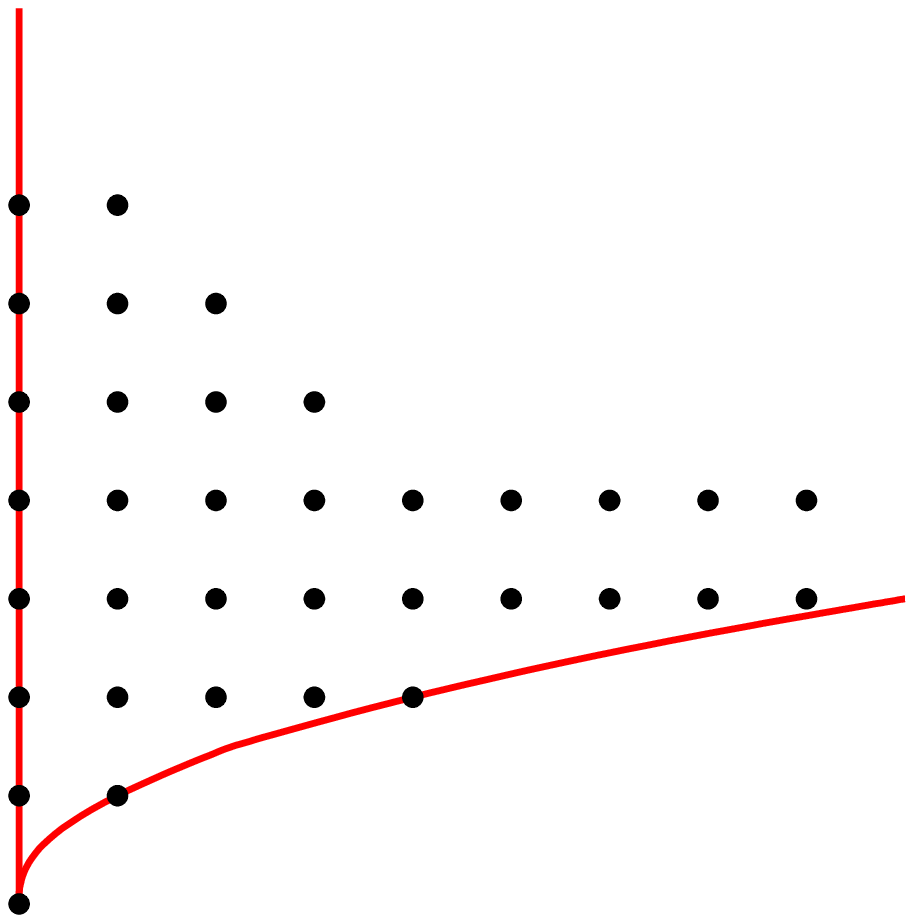, height=3cm}\\[.3cm]
\refstepcounter{figure} \label{SemiAlgCone} Figure \thefigure:
Semi-algebraic set C with its lattice points
\end{center}
\vspace*{-0.3cm}
\end{figure}

\begin{example} \label{Ex:SemiAlgCone}
The function $g$ given by
\[
g= \left(
\begin{array}{l}
\lambda_1^2\\
\lambda_1+\lambda_2\\
\end{array}
\right)
\]
defines a semi-algebraic set ${\cal
C}=\{y\in\R^2:y=g(\lambda),\lambda\in\Z^2_+\}$, see
Figure~\ref{SemiAlgCone}. In Cartesian coordinates, ${\cal C}$ can
be described by ${\cal C}=\{(x,y)\in\R^2: x-y^2\geq 0, x,y\geq
0\}$.
\end{example}

Of special interest to us will be the lattice points inside
semi-algebraic sets.

\begin{definition}
Consider a set $S\subseteq\Z^n$. Let sets $T_i\subseteq\Z^n$ be
given where each $T_i$ is described in the form
$T_i:=\{f_i(t_i):t_i\in\Z^{n_i}_+\}$ with a polynomial function
$f_i:\Z^{n_i}_+\rightarrow\Z^n$.

Then we call such a family $\{T_i\}$ an {\bf integral function
basis} of $S$, if for every $s\in S$ there exists a finite
representation, $s=\sum f_i(t_i)$, with $t_i\in\Z^{n_i}_+$ and
$f_i(t_i)\in S$.
\end{definition}

If we allowed only linear functions $f_i$ and if $S$ are the
lattice points in a rational polyhedral cone, this definition
coincides with the definition of a Hilbert basis.

If we reconsider the example with $S=\{(x,y)\in\Z^2_+:y\geq 1\}$,
we see that the following set $T_1$ defines an integral function
basis of $S$:
\[
T_1=\{f(\lambda,\mu)=(\lambda,1+\mu):\lambda,\mu\in\Z_+\}.
\]

\begin{figure}[tbh]
\begin{center}
\epsfig{file=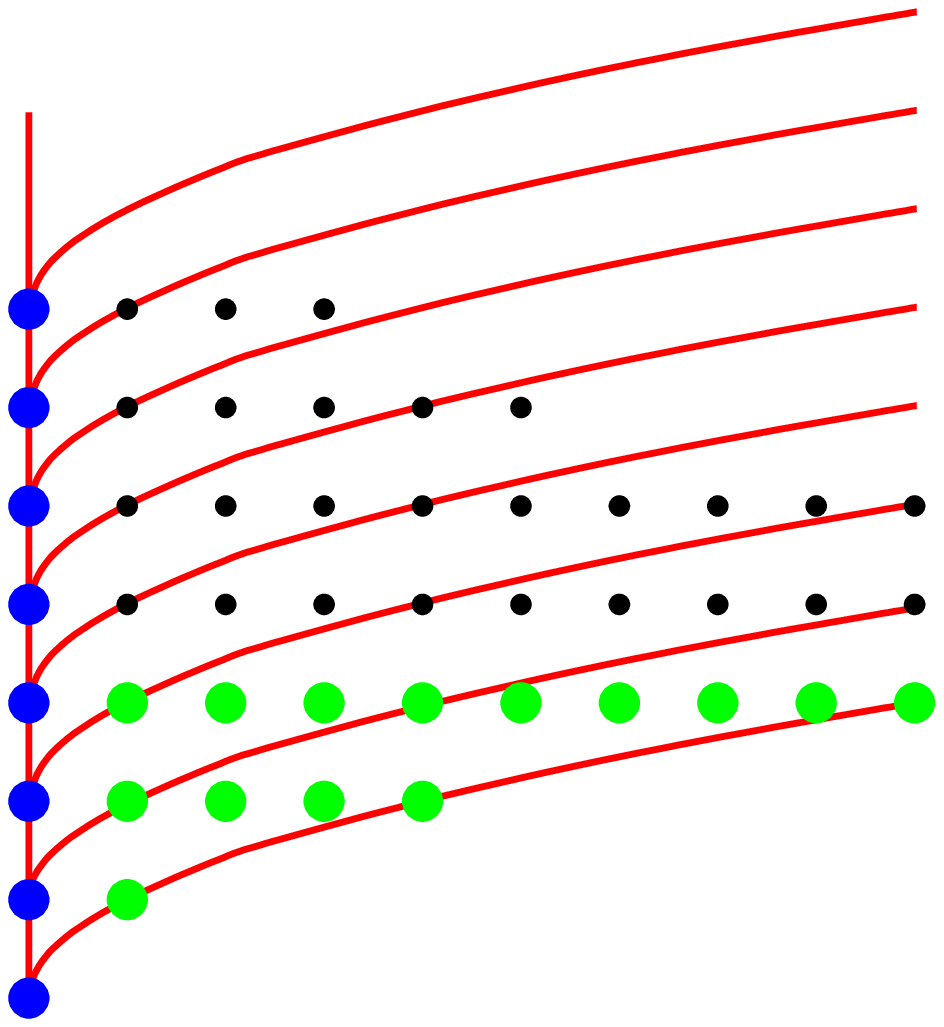, height=3cm}\\[.3cm]
\refstepcounter{figure} \label{SemiAlgCone2} Figure \thefigure:
Integral function basis of a semi-algebraic set
\end{center}
\vspace*{-0.3cm}
\end{figure}

{\bf Example \ref{Ex:SemiAlgCone}, cont.} Let us consider again
the semi-algebraic set ${\cal C}=\{(x,y)\in\R^2: x-y^2\leq 0,
x,y\geq 0\}$. An integral function basis of $S={\cal C}\cap\Z^2$
is given by $\{T_1,T_2\}$ with $T_1=\{(0,t):t\in\Z_+\}$ and
$T_2=\{((x+1)^2-s,x+1):x,s\in\Z_+\}$, where the parameters $x$ and
$s$ need to satisfy $s\leq (x+1)^2$ to guarantee
$((x+1)^2-s,x+1)\in S$, see Figure~\ref{SemiAlgCone2}. The only
lattice point in ${\cal C}$ that cannot written as a sum of a
lattice point in $T_1$ and a lattice point in $T_2$ is the origin.
This special point, however, can already be represented by $T_1$
alone.

\begin{figure}[tbh]
\begin{center}
\epsfig{file=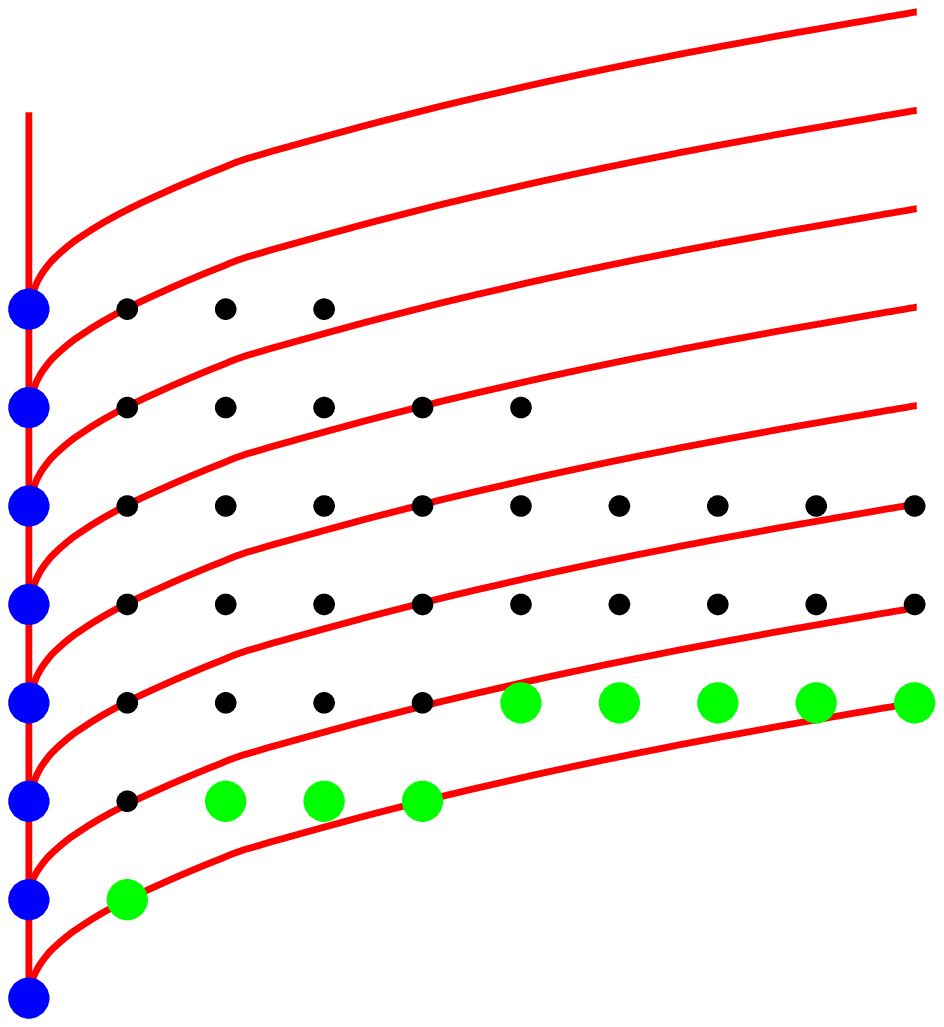, height=3cm}\\[.3cm]
\refstepcounter{figure} \label{SemiAlgCone3} Figure \thefigure:
Integral function basis of a semi-algebraic cone
\end{center}
\vspace*{-0.3cm}
\end{figure}

Note that the constraints that are needed to encode the condition
$(x+1)^2-s,x+1)\in S$ have the same maximal degree as the original
constraints. On the other hand, we can also this condition by
$s\leq 2x+1$, see Figure~\ref{SemiAlgCone3}. The latter
representation should be preferred, since this description of the
set $S$ using the integral function basis involves only linear
constraints in contrast to the quadratic constraint in the
description above. \eoproof

In the following, we outline a fundamental application of integral
function bases for nonlinear integer optimization problems. It
turns out that one can derive an optimality criterion for a linear
integer program with a polynomial objective function. This
criterion is a natural generalization of what Graver proved in the
fully linear (integer) setting \cite{Graver:75}.

\begin{theorem} \label{Theorem: Application to integer
optimization} Let $p(z)$ be any polynomial in
$\Z[z_1,\ldots,z_n]$, and let $A\in\Z^{d\times n}$ and $b\in\Z^d$.
Consider the (possibly nonlinear) integer optimization problem
\[
\max\{p(z):Az=b,z\geq 0, z\in\Z^n\}.
\]
Let $W^i\in\Z^{n\times w_i}$, $i=1,\ldots 2^n$, denote the extreme
rays of the cones
\[
\cone(W^i)=\{x\in\R^n\cap\Orthant_i:\;Ax=0\},
\]
where $\Orthant_1,\ldots,\Orthant_{2^n}$ denote the $2^n$ orthants
of $\R^n$. Thus, every point in this cone can be written as a
linear combination $z=W^i\lambda$, $\lambda_i\geq 0$.

Assume that $z_0$ is a feasible integer solution to $Az=b$, $z\geq
0$. For each $i=1,\ldots,2^n$, define the following vector of
nonlinear functions,
\[
g^i(\lambda)= \left(
\begin{array}{c}
g^i_1(\lambda)\\
\vdots\\
g^i_n(\lambda)\\
\bar{q}^i(\lambda)\\
\end{array}
\right)= \left(
\begin{array}{c}
W^i\lambda\\
p(z_0+W^i\lambda)-p(z_0)\\
\end{array}
\right).
\]
Let $\{T^i_1,\ldots,T^i_{k_i}\}$, with
$T^i_j:=\{f^i_j(t_j):t_j\in\R^{n_{i,j}}_+\}$, be an integral
function basis for the (integer points in the) semi-algebraic set
${\cal C}^i= \{y\in\R^{n+1}:y=g^i(\lambda),\lambda\geq 0\}$.
Define by $S_{i,j}$ the semi-algebraic set that encodes the
conditions $[f^i_j(t_j)]_{1,\ldots,n}\in\cone(W^i)$. (Herein,
$[f^i_j(t_j)]_{1,\ldots,n}$ shall denote the vector of the first
$n$ components of $f^i_j(t_j)$.)

Then $z_0$ is optimal if and only if for every $i=1,\ldots,2^n$,
the following condition holds:
\[
[f^i_j(t_j)]_{n+1}\leq 0 \text{ for all } t_j\in\Z^{n_{i,j}}_+\cap
S_{i,j} \text{ with } [f^i_j(t_j)]_k\geq -[z_0]_k, \text{ for all
} k=1,\ldots,n.
\]
\end{theorem}

\boproof Assume that there is a better feasible solution $z_1$
that has an objective value $p(z_1)>p(z_0)$. Consider the
difference vector $v:=z_1-z_0$, which lies in one of the $2^n$
orthants $\Orthant_i$ of $\R^n$. Therefore, we are looking for
$v\in\Orthant_i$ with $Av=0$, $z_0+v\geq 0$, and
$p(z_0+v)-p(z_0)>0$. Clearly, the set $\{z\in\Orthant_i:Az=0\}$
forms a pointed rational cone, generated by the columns of $W^i$.
Thus, $v=W^i\lambda$ for some $\lambda\in\R^{w_i}_+$ and hence
\[
\left(
\begin{array}{c}
W^i\lambda\\
p(z_0+W^i\lambda)-p(z_0)\\
\end{array}
\right)
\]
is an integer point in the semi-algebraic set ${\cal C}^i=
\{y\in\R^{n+1}:y=g^i(\lambda),\lambda\geq 0\}$. Using the integral
function basis of this set, there is a representation
\[
\left(
\begin{array}{c}
W^i\lambda\\
p(z_0+W^i\lambda)-p(z_0)\\
\end{array}
\right)= \sum\limits_{j\in I^i} f^i_j(t_j)= \sum\limits_{j\in I^i}
\left(
\begin{array}{r}
[f^i_j(t_j)]_{1,\ldots,n}\\
{[f^i_j(t_j)]}_{n+1}\\
\end{array}
\right)
\]
with $t^i_j\in S_{i,j}\cap\Z^{n_{i,j}}$ and $f^i_j(t_j)\in {\cal
C}^i$.

As $p(z_0+W^i\lambda)-p(z_0)>0$, there must be some $j\in I^i$
with $[f^i_{j}(t_j)]_{n+1}>0$. We claim that the first $n$
components of $f^i_j(t_j)\in {\cal C}^i$ form an improving integer
vector for $z_0$, possibly different from the vector $v$ that we
decomposed.

As $[f^i_{j}(t_j)]_{n+1}>0$, the only thing left to show is that
the components $[f^i_j(t_j)]_1,\ldots,[f^i_j(t_j)]_n$ of
$f^i_j(t_j)$ lie above the lower bounds, i.e., $[f^i_j(t_j)]_k\geq
-[z_0]_k$ for $k=1,\ldots,n$. But this can be seen as follows. By
construction, $\cone(W^i)\in\Orthant^i$, implying
$[f^i_j(t_j)]_{1,\ldots,n}\in\Orthant^i$ for all $j$. Thus, the
components of $z_0+[f^i_j(t_j)]_{1,\ldots,n}$ lie between the
components of $z_0$ and of $z_0+v=z_1$, and are therefore
nonnegative.

The converse direction is obviously true. \eoproof

Clearly, one would wish that searching for an improving vector in
each of the $T_i$ is simpler than searching for an improving
vector in ${\cal C}$.

The set $T_1=\{\lambda-\mu:\lambda,\mu\in\Z^n_+\}$ always forms an
integral function basis for any set $S\subseteq\Z^n$, where $2n$
parameters are needed to describe $T_1$. The following theorem
bounds the number of parameters needed in the $T_i$ and thus gives
a sufficient condition (together with a construction) of when an
integral function basis with less parameters in the description of
each $T_i$ exists.

\begin{theorem} \label{Theorem: Every set S has in IFB}
Let $S\subseteq\Z^n$, $v_1,\ldots,v_k\in\Z^n$ and let
$C=\cone(v_1,\ldots,v_k)$ be a rational polyhedral cone with
$S\subseteq C$. Then $S$ has an integral function basis in which
the appearing sets $T_i$ involve at most $k+1$ parameters.
\end{theorem}

\boproof First observe that $S\subseteq C$ implies
$\cone(S)\subseteq C$ and therefore
$\cone(S\cup\{v_1,\ldots,v_k\})=C$. Thus, by Theorem \ref{General
Finiteness Theorem}, there is a finite integral basis
$\{h_1,\ldots,h_s,v_1,\ldots,v_k\}$ for the set
$S\cup\{v_1,\ldots,v_k\}$. If we set in addition $h_0=0$, we can
see from the proof of Theorem \ref{General Finiteness Theorem},
that every point $v\in S$ can be written as
$v=h_i+\sum_{j=1}^k\lambda_{ij} v_j$ for some $i\in\{0,\ldots,s\}$
and for some nonnegative integers $\lambda_{ij}$. This last
condition in fact states that the sets $T_i=\{h_i+\sum_{j=1}^k
\lambda_{ij} v_j:\lambda_{i}\in\Z^k_+\}\cup\{0\}$,
$i=0,1,\ldots,s$, form an integral function basis for $S$.
\eoproof

\begin{remark}
It should be noted that we may strengthen the above theorem if
some or all of the cone generators $v_j$ lie in $S$. If
$v_{j_0}\in S$, then each set $T_i=\{h_i+\sum_{j=1}^k \lambda_{ij}
v_j:\lambda_{i}\in\Z^k_+\}\cup\{0\}$ can in fact be decomposed
into the sum of
$T'_i=\{h_i+\sum_{j\in\{1,\ldots,k\}\setminus\{j_0\}} \lambda_{ij}
v_j:\lambda_i\in \Z^{k-1}_+\}\cup\{0\}$ and
$T''_i=\{\lambda_{ij_0} v_{j_0}:\lambda_{ij_0}\in\Z_+\}$.

Iterating this process for all cone generators $v_j$ that lie in
$S$ gives a new integral function basis for $S$ with fewer
parameters appearing in the description of the sets $T_i$. In
fact, if all $v_j$ lie in $S$, that is if $C=\cone(S)$, the
integral function basis for $S$ simplifies to sets $T_i$ that all
contain nonnegative integer multiples of a single lattice point of
$S$. Thus, we have recovered the statement of Theorem \ref{General
Finiteness Theorem}: the existence of a finite integral basis if
$\cone(S)$ is rational.
\end{remark}

The following example demonstrates that splitting the set $S$ into
finitely many subsets may also decrease the maximum number of
parameters needed in the description of the $T_i$'s.

\begin{example}
Consider the set $S=\{(x,y)\in\Z^2: -y^2\leq x\leq y^2\}$. As this
set is contained in the rational cone spanned by $e_1$ and $-e_1$,
we conclude by Theorem \ref{Theorem: Every set S has in IFB} that
$S$ has an integral function basis, in which each $T_i$ is
described by at most $2+1=3$ parameters.

However, if we split the set $S$ as
\[
S=S'+S''=\{(x,y)\in\Z^2: -y^2\leq x\leq 0\}\cup \{(x,y)\in\Z^2:
0\leq x\leq y^2\},
\]
we see that $S'\subseteq\cone(-e_1,e_2)$ and
$S''\subseteq\cone(e_1,e_2)$. Since the cone generator $e_2$ is an
element of the sets $S'$ and $S''$, respectively, we find integral
function bases for $S'$ and $S''$, in which each $T_i$ is
described by at most $1+1=2$ parameters. Putting both together, we
arrive at an integral function basis for $S$ with the same
property. \eoproof
\end{example}

As we have seen above, every set of lattice points in $\Z^n$
admits a representation via an integral function basis. Even under
the assumption that we have found a nice integral function basis
for a particular problem instance, that is, one that has only few
parameters in the description of the $T_i$, we are faced with a
new problem to be solved.

Suppose we want to maximize a (polynomial) function $p(x)$ over
the lattice points in a semi-algebraic set ${\cal C}$. Knowing an
integral function basis $\{T_i:i\in I\}$, we can use the
representation $x=\sum_{i\in I} f_i(t_i)$, $t_i\in\Z^{n_i}$ for
all $x\in {\cal C}\cap\Z^n$ to rewrite the problem as
\[
\max\left\{p\left(\sum_{i\in I} f_i(t_i)\right):\sum_{i\in I}
f_i(t_i)\in {\cal C}\cap\Z^n, t_i\in\Z^{n_i}_+\right\}.
\]
While the condition $\sum_{i\in I} f_i(t_i)\in {\cal C}\cap\Z^n$
often follows immediately from $f_i(t_i)\in {\cal C}\cap\Z^n$ for
all $i\in I$, these latter conditions involve descriptions by
polynomials of the same degree as in the description of ${\cal C}$
and are thus still hard to deal with. Finding $t_i\in\Z^{n_i}_+$
with $f(t_i)\in {\cal C}\cap\Z^n$ even only for a single $i$ (as
needed in Theorem \ref{Theorem: Application to integer
optimization}) is as hard as finding a point in ${\cal C}$, at
least from a complexity point of view.

Thus, an integral function basis with the additional property that
$T_i\subseteq {\cal C}\cap\Z^n$ for all $i\in I$ would be
desirable. Then $f_i(t_i)\in {\cal C}\cap\Z^n$ for all $i\in I$
and for all $t_i\in\Z^{n_i}$ would hold automatically. For this,
of course, a nonlinear description for the $T_i$ is needed, in
contrast to the rather nice and simple description guaranteed to
exist by Theorem \ref{Theorem: Every set S has in IFB}.

In the following, we relax the condition $f_i(t_i)\in {\cal
C}\cap\Z^n$ and allow a {\it correction term} that may lie
outside, but which is bounded by polynomials of strictly smaller
degree than the given polynomials.

\begin{theorem}
For every semi-algebraic set
\[
{\cal C}:=\{x\in\R^n:\exists y\geq 0 \text{ with } x=g(y)\},
g:\R^d\rightarrow\R^n,
\]
there exists a set of functions
\[
\{g_l,g_u:\R^n\rightarrow\R^n\}
\]
with $\maxdeg(g_l),\maxdeg(g_u)<\maxdeg(g)$ and such that for
every point $x\in {\cal C}\cap\Z^n$ there exists a
$\lambda\in\Z^d_+$ and a point $v_x\in\Z^n$ with
$x=g(\lambda)+v_x$ and with $g_l(\lambda)\leq v_x\leq
g_u(\lambda)$.
\end{theorem}

\boproof Choose any $x\in S:={\cal C}\cap\Z^n$. Then $x=g(y)$ for
some $y\in\R_+^n$. Now define $\lambda:=\lfloor y\rfloor$
component-wise and let $v_x=x-g(\lambda)$ and $h=y-\lambda$. We
will now construct functions $g_{l}:\R^n\rightarrow\R^n$ and
$g_{u}:\R^n\rightarrow\R^n$, with the desired properties.

Let $D=\maxdeg(g)$. By multivariate Taylor expansion, we get for
$j=1,\ldots,n$:
\[
x^{(j)}=g^{(j)}(\lambda+h)=g^{(j)}(\lambda)+\sum_{i=1}^D
\frac{1}{i!}\cdot
\sum_{\begin{array}{c}\alpha\in\Z_+^n:\\\|\alpha\|_1=i\\\end{array}}
\frac{{\rm d} g^{(j)}(\lambda)} {{\rm d} x^\alpha}\cdot h^\alpha.
\]
Therefore,
\[
v_x^{(j)}=x^{(j)}-g^{(j)}(\lambda)=
\sum_{i=1}^D\frac{1}{i!}\cdot
\sum_{\begin{array}{c}\alpha\in\Z_+^n:\\\|\alpha\|_1=i\\\end{array}}
\frac{{\rm d} g^{(j)}(\lambda)} {{\rm d} x^\alpha}\cdot h^\alpha.
\]
Note that
$\maxdeg\left(\frac{{\rm d} g^{(j)}} {{\rm d} x^\alpha}\right)<
\maxdeg(g)$ and that $0\leq h<1$ by construction.

This sum is a polynomial in $\lambda$ and $h$, that is, it is a sum of
terms $c_{\alpha,\beta} \lambda^\alpha h^\beta$. Since all
$\lambda\geq 0$ we can use $0\leq h_i<1$, for all $i$, to bound the
expression $c_{\alpha,\beta} \lambda^\alpha h^\beta$ by
\[
0\leq c_{\alpha,\beta} \lambda^\alpha h^\beta<
c_{\alpha,\beta} \lambda^\alpha
\]
if $c_{\alpha,\beta}>0$ and by
\[
c_{\alpha,\beta} \lambda^\alpha < c_{\alpha,\beta} \lambda^\alpha
h^\beta \leq 0
\]
if $c_{\alpha,\beta}<0$.
Putting now
\[
g_{l}^{(j)}(\lambda):=\sum_{\alpha,\beta:c_{\alpha,\beta}<0}
c_{\alpha,\beta} \lambda^\alpha\;\;\;\text{ and }\;\;\;
g_{u}^{(j)}(\lambda):=\sum_{\alpha,\beta:c_{\alpha,\beta}>0} c_{\alpha,\beta} \lambda^\alpha
\]
we have
\[
g_{l}^{(j)}(\lambda)\leq v_x^{(j)}\leq g_{u}^{(j)}(\lambda)
\]
by construction. Moreover, again by construction, the degree of
$g_{l}^{(j)}$ and of $g_{u}^{(j)}$ is strictly less than the degree of
$g^{(j)}$.
\eoproof

The above theorem tells us that the error term $v_x$ can be
bounded by polynomials of strictly smaller maximal degree than
that of $g(\lambda)$. As the following example shows, the degree
of $g_{l}^{(j)}$ and of $g_{u}^{(j)}$ can in fact be much smaller
than that of $g^{(j)}$.

\begin{example} Let us consider again the semi-algebraic set given by
\[
g(y)=\left(
\begin{array}{c}
y_1\\
y_1^k+y_2\\
\end{array}
\right).
\]
As can be easily checked, each integral point $v$ in this
semi-algebraic set can be written as $v=g(\lambda)$ for
$\lambda\in\Z^2_+$, showing that the correction term $v_x$ is $0$
in this case. \eoproof
\end{example}

This leads us immediately to the questions of when is $v_x=0$ or
of when is $v_x\in {\cal C}\cap\Z^n$? In both cases, of course,
$T=\{g(\lambda):\lambda\in\Z^d_+\}$ would be an integral function
basis for ${\cal C}$ with our desired property $T\subseteq {\cal
C}\cap\Z^n$.

We believe that research in this direction will make it possible
to design novel algorithms for polynomial integer programming
based on reformulation techniques.

\bibliographystyle{plain}
\bibliography{SemiAlgOpt}

\begin{thebibliography}{10}

\bibitem{Bertsimas+Weismantel}
D.~Bertsimas and R.~Weismantel.
Optimization over Integers. Manuscript,
2004, Dynamic Ideas, Belmont Mass, to appear 2005.

\bibitem{Giles+Pulleyblank:79}
F.~R.~Giles and W.~R.~Pulleyblank. Total dual integrality and
integer polyhedra. Linear Algebra and its Applications 25 (1979),
191--196.

\bibitem{Graver:75}
J.~E.~Graver. On the foundation of linear and integer programming
{I}. \textit{Mathematical Programming} \textbf{9} (1975), 207-226.

\bibitem{Haus+Koeppe+Weismantel}
U.~U.~Haus, M.~K{\"o}ppe, and R.~Weismantel.
The integral basis method for integer programming.
\textit{Mathematical Methods of Operations Research}, \textbf{53}
(2001).

\bibitem{lenstra}
Lenstra, H.W.: Integer Programming with a fixed number of
variables. {Mathematics of Operations Research}, 8, (1983)
538--548.

\bibitem{Schrijver:86}
A.~Schrijver. \textit{Theory of Linear and Integer Programming}.
Wiley, Chichester, 1986.

\end{thebibliography}

\end{document}